\documentclass[reqno]{amsart}
\usepackage[pagebackref]{hyperref}
\usepackage{amssymb}
\usepackage{amscd}
\input{diagrams}

\newcommand{\bbeta}{\overline{\beta}}
\newcommand{\hbeta}{\widehat{\beta}}
\newcommand{\bgamma}{\overline{\gamma}}
\newcommand{\hgamma}{\widehat{\gamma}}
\newcommand{\hg}{\widehat{g}}

\newcommand{\lra}{\longrightarrow}
\newcommand{\too}{\longmapsto}
\newcommand{\id}{\operatorname{id}}
\newcommand{\im}{\operatorname{im}}
\newcommand{\cok}{\operatorname{coker}}
\newcommand{\coker}{\operatorname{coker}}
\newtheorem{thm}{Theorem}

\newtheorem{cor}[thm]{Corollary}

\begin{document}
\title{The Snake Lemma}
\author{Franz Lemmermeyer}
\address{M\"orikeweg 1, 73489 Jagstzell}
\email{hb3@ix.urz.uni-heidelberg.de}

\maketitle

The snake lemma in abelian categories is a simple and  
very useful result; in the following, we will present 
a version of the snake lemma that contains the usual 
formulation as a special case. 

\begin{thm}[Snake Lemma]\label{SL}
Assume that
$$\begin{CD}
   A   @>>f>     B  @>>g>  C \\     
 @VV{\alpha}V  @VV{\beta}V  @VV{\gamma}V   \\
   A'  @>>{f'}>  B' @>>g'> C' \end{CD}$$
is a commutative diagram of abelian groups with exact rows. Then  
there exists a homomorphism $\delta: \ker \gamma \, \cap \, \im g 
\lra A' / (\im \alpha + \ker f')$ 
such that the following sequence is exact:
$$\begin{CD}
0 @>>> \ker f @>>> \ker f' \circ \alpha @>>> \ker \beta @>>> 
               \ker \gamma \, \cap \, \im g \\ @. @. @. @. @V{\delta}VV \\
0  @<<< \cok g' @<<<  \cok \gamma \circ g @<<< \cok \beta 
   @<<< A'/(\im \alpha + \ker f') \end{CD} $$
If $f'$ is injective, then $\ker f' \circ \alpha = \ker \alpha$ 
and $A'/(\im \alpha + \ker f') = \cok \alpha$; if $g$ is 
surjective, then $\cok \gamma \circ g = \cok \gamma$ and 
$\ker \gamma \, \cap \, \im g = \ker \gamma$. Thus if $f'$ is 
injective and $g$ is surjective, then we get the following 
exact sequence:
$$\begin{CD}
    0 @>>> \ker f @>>> \ker \alpha @>>> \ker \beta @>>> \ker \gamma \\ @.
      @. @. @. @V{\delta}VV \\
    0  @<<< \cok g' @<<<  \cok \gamma  @<<< \cok \beta @<<< \cok \alpha
  \end{CD} $$
\end{thm}

The proof of the standard version of the snake lemma goes through.

\begin{cor}[Ring Lemma] \label{C1}
Let $\alpha: A \lra B$ and $\beta:B \lra C$ be homomorphisms; then
there is an exact sequence
$$\begin{CD}
  0  @>>> \ker \alpha @>>> \ker (\beta \circ \alpha) @>>> 
	\ker \beta \\
   @. @. @. @VVV \\
  0  @<<<\coker \beta  @<<<  \coker (\beta \circ \alpha)  @<<< 
	\coker \alpha 
\end{CD}$$
\end{cor}

\begin{proof}
Apply the snake lemma to the diagram
$$ \begin{CD}
 @. A  @>{\alpha}>> B @>>> \coker \alpha @>>> 0 \\
 @.  @VV{\beta \circ \alpha}V  @VV{\beta}V  @VVV  @.   \\
 0 @>>> C @>{\id}>> C  @>>> 0   @. 
\end{CD} $$
\end{proof}

Bass 
has observed that the $6$-term exact sequence fits into the 
following exact and commutative diagram (the exact ring):

\begin{diagram}[width=2.5em,height=1.6em] 
     & &   & \ker \beta &       \rTo & \coker \alpha &   \\ 
     & & \ruTo(2,4) &            & \rdInto(1,2) \ruOnto(1,2) & 
     & \rdTo(2,4) & &   \\
     & & &            &       B    & &    \\  
     &      &            &       \ruTo(1,2) & & \rdTo(1,2) &  &  &  & \\  
  & \ker \beta \circ \alpha & \rInto & A &      \rTo 
     &       C & \rOnto & \coker \beta \circ \alpha  \\
     &   & \luTo(1,2) \ruInto (1,2) &  & &  & \rdOnto(1,2) \ldTo (1,2) \\  
     &   &          \ker \alpha & \lTo  & 0 & \lTo & \coker \beta  \\
\end{diagram}

\begin{cor}[$4$-Lemma]
Assume that the diagram
$$ \begin{CD} A  @>f>>  B @>g>> C @>h>> D \\
   @V{\alpha}VV @V{\beta}VV @V{\gamma}VV @V{\delta}VV \\
   A' @>f'>> B' @>g'>> C' @>h'>> D' \end{CD} $$
is commutative with exact rows. If $\alpha$ is surjective and
if $\delta$ is injective, then we have the following exact sequences

\begin{equation} \label{ES1} \begin{CD}  
 0 @>>> \ker \beta \cap \ker g @>>> \ker \beta @>{g^*}>> \ker \gamma @>>> 0,
\end{CD} \end{equation}
\begin{equation}\label{ES2} \begin{CD} 
 0 @>>> \cok \beta @>{g_*'}>> \cok \gamma @>>> C'/(\im \gamma + \im g') @>>> 0.
 \end{CD} \end{equation}

In particular, we have $\ker \gamma = g(\ker \beta)$ and 
$\im \beta = {g'}^{-1}(\im \gamma)$.
\end{cor}

\begin{proof}
Apply the snake lemma to the diagram consisting of the second 
and the third square; observing that $\ker \delta = 0$ provides 
us with the exact sequence
\begin{equation}\label{ESR} \begin{CD} 
   0 @>>> \ker g @>{\iota}>> \ker g' \circ \beta @>{\hg}>> 
          \ker \gamma @>>> 0. \end{CD} \end{equation}
Next, $\beta$ induces a map $\hbeta: \ker g' \circ \beta \lra \ker g'$.
Using the fact that $\alpha$ is surjective we find
$\beta(\ker g) = \beta(\im f) = \im \beta \circ f  = 
 \im f' \circ \alpha = \im f' = \ker g'$, hence 
$\cok \hbeta \circ \iota = 0$; applying the ring lemma to $\iota$
and $\hbeta$ and observing that $\ker \iota = 0$ we get the exact 
sequence
$$ \begin{CD} 
   0 @>>> \ker \hbeta \circ \iota @>>> \ker \hbeta 
     @>>> \cok \iota @>>> 0. \end{CD} $$
This gives (\ref{ES1}) since 
$\ker \hbeta \circ \iota = \ker \beta \cap \ker g$,
$\ker \hbeta = \ker \beta$, and finally 
$\cok \iota = \ker \gamma$ by (\ref{ESR}).

The proof of the dual sequence (\ref{ES2}) is similar.
\end{proof}
\end{document}